\documentclass[letterpaper, 10 pt, conference]{ieeeconf}  

\usepackage[utf8]{inputenc}
\usepackage[english]{babel}
\usepackage{graphicx}
\usepackage{amsmath,amsfonts}
\usepackage{epstopdf}
\usepackage{graphicx}
\graphicspath{{images/}{../images/}}
\usepackage{subfiles}  
\usepackage{blindtext} 

\newtheorem{definition}{Definition}
\newtheorem{lemma}{Lemma}
\newtheorem{theorem}{Theorem} 
\newtheorem{proposition}{Proposition}
\newtheorem{assumption}{Assumption}
\newtheorem{remark}{Remark}

\usepackage[usenames, dvipsnames]{color}
\usepackage{ifthen}
\newboolean{showcomments}
\setboolean{showcomments}{false}   
\newcommand{\enrique}[1]{  \ifthenelse{\boolean{showcomments}}
{\textcolor{red}{(Enrique says:  #1)}}{}}
\newcommand{\addcite}[0]{\ifthenelse{\boolean{showcomments}}
{\textcolor{Purple}{~(add cite(s)) } }{}}
\newcommand{\addcites}[0]{\ifthenelse{\boolean{showcomments}}
{\textcolor{Purple}{~(add cite(s)) } }{}}

\newcommand{\tomove}[1]{\ifthenelse{\boolean{showcomments}}{#1}{}}
\usepackage{setspace} 

\IEEEoverridecommandlockouts                               
\title{\LARGE \bf An integral quadratic constraint framework for real-time steady-state optimization of linear time-invariant systems*}

\author{Zachary E. Nelson and Enrique Mallada
\thanks{*This work was supported by the Army Research Office contract W911NF-17-1-0092, NSF grants (CNS 1544771, EPCN 1711188), and Johns Hopkins WSE startup funds.} 
\thanks{Zachary E. Nelson and Enrique Mallada are with the Department of Electrical and Computer Engineering, The Johns Hopkins University, 3400 N. Charles Street, Baltimore, MD 21218, emails:
	{\tt\small \{znelson2,mallada\}@jhu.edu}}
} 
 
\begin{document}
 
\maketitle
\thispagestyle{empty}
\pagestyle{empty}

\begin{abstract} 
Achieving optimal steady-state performance in real-time is an increasingly  necessary requirement of many critical infrastructure systems. In pursuit of this goal, this paper builds a systematic design framework of feedback controllers for Linear Time-Invariant (LTI) systems that continuously track the optimal solution of some predefined optimization problem. The proposed solution can be logically divided into three components. The first component estimates the system state from the output measurements. The second component uses the estimated state and computes a drift direction based on an optimization algorithm. The third component computes an input to the LTI system that aims to drive the system toward the optimal steady-state.

We analyze the equilibrium characteristics of the closed-loop system and provide conditions for optimality and stability. Our analysis shows that the proposed solution guarantees optimal steady-state performance, even in the presence of constant disturbances. Furthermore, by leveraging recent results on the analysis of optimization algorithms using integral quadratic constraints (IQCs), the proposed framework is able to translate input-output properties of our optimization component into sufficient conditions, based on linear matrix inequalities (LMIs), for global exponential asymptotic stability of the closed loop system. We illustrate the versatility of our framework using several examples. 

\end{abstract}

\section{Introduction}\label{sec:introduction}

Infrastructure systems are the foundation of our modern society. The Internet, power grids, and transportation networks are just some examples of the several critical systems that our current lifestyle relies on. Due to their large scale, the operational and fault-associated costs that these systems incur are both in the range of hundreds of millions of dollars to several billion dollars \cite{Campbell2012}. Therefore, operators are continuously faced with the conflicting tasks of operating these systems as efficiently as possible and guaranteeing certain levels of security or robustness.

Traditionally, this balancing between efficiency and security is achieved by separating tasks across different time-scales. Efficiency goals are achieved using optimization algorithms running at a slow time-scale, and stability/robustness goals are achieved using fast time-scale controllers. For example, in power systems, generators are optimally scheduled by solving an (economic dispatch) optimization problem at a slow time-scale (every 5/15 minutes, hour, or day) \cite{Bejestani2014}, but at the fast time-scale the scheduling uses controllers based on frequency measurements \cite{Wood2012} that are focused on preserving the system stability \cite{Ibraheem2005}, not efficiency.

Unfortunately, the state of flux that these infrastructure systems currently experience due to the growing population, deployment of sensing and communication technologies, and sustainability trends, is pushing its operation towards their limits and, in this way, rendering this approach obsolete. Operating at maximum capacity does not leave room for the inefficiencies incurred by the timescale separation. Moreover, the limited coordination capabilities that today's controllers provide, when compared with optimization algorithms, does not allow the system to quickly react to unprescribed events. Motivated by this problem, this paper aims to remove the time-scale separation by building controllers that can simultaneously achieve steady-state optimality while preserving the system stability.

More precisely,  this paper proposes a systematic design framework for feedback controllers that, given a LTI system and an unconstrained optimization problem,  generates a family of nonlinear controllers that seek to drive the system towards the optimal solution of the optimization problem.
Our equilibrium analysis connects notions of controllability and observability of the LTI system with a criterion for steady-state optimality of the closed loop equilibrium. Furthermore, we leverage recent analyses of optimization algorithms using Integral Quadratic Constraints (IQCs)~\cite{Lessard2016,Fazlyab2017,Hu2017}  to provide sufficient conditions based on Linear Matrix Inequalities~\cite{Boyd1994} that guarantee global exponential asymptotic stability. 
The derived LMIs provide an explicit bound on the rate of convergence and allow us to design an algorithm that computes the maximum rate of convergence.

Our controllers have two main distinctive features. Firstly, they can be functionally separated into three components/modules: (i) an Estimator, that aims to estimate the system's state from the output; (ii) an Optimizer, that uses the estimated state to compute the drift direction necessary to achieve optimality or outputs zero when the estimated state is optimal; and (iii) a Driver (PI controller) that generates the necessary input to drive the system toward the optimal solution. Secondly, the Optimizer module can be implemented using one of many optimization algorithms, leading to a family of optimization-based nonlinear controllers. The only required conditions are that (i) in steady state the output of the optimizer is zero if and only if its input (the estimated state) is the optimal solution of the optimization problem, and (ii) there exists an Integral Quadratic Constraint that captures the input-output relationship of the optimizer.

\noindent
\emph{Related Work:}
Optimization-based control design for achieving optimal steady-state performance has a been a popular subject of research for more than three decades. It has been used in communication networks to reverse engineer TCP/IP congestion control protocols~\cite{Kelly98, Low99} and provide a design framework for novel congestion control algorithms~\cite{Wei2006}, distributed multi-path routing~\cite{Rieck2005,Mallada2008}, and admission control~\cite{Ferragut2008}, and access control in wireless networks~\cite{Chen2006}. In the context of power systems and micro-grids, optimization-based control design has been used for the design of distributed controllers that can achieve efficient supply-demand balance~\cite{Zhao:2014bp}, frequency restoration~\cite{ml2014ifacwc,mzl2017tac}, congestion management~\cite{zmbl2016pscc}, and economic steady-state optimality~\cite{Dorfler:2014uu,Cherukuri2015,zmd2015acc,Li2016}. Some of these approaches have been further extended for more general settings such as~\cite{Jokic2009} and \cite{Zhang2015}. In general, these solutions either require that the dynamical system to be optimized has a specific structure, such as being passive~\cite{zmbl2016pscc,Dorfler:2014uu}, having primal-dual dynamics~\cite{mzl2017tac,Li2016,Zhang2015}, or having direct access to (a subset of) the system state~\cite{Jokic2009}. 

More recently, real-time optimization algorithms have been proposed as a mean to mitigate the large fluctuations that renewable energy introduce in power networks. The solutions fall within two categories depending on whether the system dynamics are considered as perturbations of the optimization algorithms~\cite{DallAnese2016Optimal,DallAnese2016Photovoltaic}, or the system is modeled as a set of nonlinear algebraic constraints with slowly time varying parameters~\cite{Dorfler2016, Low2017}. Our work distinguishes from these works by explicitly modeling the system dynamics and simultaneously guaranteeing stability of the dynamical system and convergence to the optimal solution. Notably, while our framework today does not include optimization constraints or nonlinearities in the system dynamics, extending our framework to incorporate these features is the subject of our current research.

\noindent
\emph{Paper Organization:}
The organization of the paper is as follows. Section \ref{sec:preliminaries} gives the reader preliminary tools that are necessary for the later analysis. Section \ref{sec:problem setup} sets up the problem and discusses some of the challenges. Section \ref{sec:control design} proposes a design framework of controllers that addresses the challenges. Section \ref{sec:analysis} shows the systematic procedure for analyzing steady-state optimality and stability. Section \ref{sec:numerical lllustrations} considers multiple numerical examples to illustrate the practicality of this approach. Lastly, Section \ref{sec:conclusions} summarizes the major points of the paper and suggests future work.


\section{Preliminaries}\label{sec:preliminaries}

\subsection{Notation}
The following notation will be used throughout the remainder of the paper. The $n \times n$ identity matrix is denoted as $I_n$. The $m \times n$ zero matrix is denoted as $0_{m \times n}$. The zero vector with length $n$ is denoted as $0_n$. The subscripts are removed when the dimensions are implied by context. A positive (semi) definite matrix $P \in \mathbb{R}^{n \times n}$ is denoted as $P \succ 0 \ (\succeq 0)$. All norms $||\cdot||: \mathbb{R}^n \rightarrow \mathbb{R}$ are the standard $\ell2$-norm. The Kronecker product of two matrices is denoted by the symbol $\otimes$.

\subsection{Integral Quadratic Constraints}

Given a nonlinear mapping $\phi: p\mapsto q$, with $p,q\in\mathbb{R}^n$, and an input-output reference $(p_*, \phi(p_*))\in \mathbb{R}^n\times \mathbb{R}^n$, we consider the following class of IQCs:
\begin{definition}[Pointwise IQC]~\label{PointwiseIQC}
	The mapping $\phi$ is said to satisfy the \emph{pointwise IQC} defined by $(Q_\phi,p_*,\phi(p_*))$ if 
	\[
	\begin{bmatrix} p-p_* \\ \phi(p)-\phi(p_*) \end{bmatrix}^T 
	Q_\phi 
	\begin{bmatrix} p-p_* \\ \phi(p)-\phi(p_*) \end{bmatrix} \geq 0 
	\]
	holds for all $(p,p_*)\in \mathbb{R}^n \times \mathbb{R}^n$, where $Q_\phi^T= Q_\phi \in \mathbb{R}^{2n \times 2n}$ is an indefinite matrix. 
\end{definition} 

Next, we discuss two particular choices of the nonlinear map $\phi$ that are commonly used in optimization algorithms.

\vspace{1ex}
\noindent
\emph{Gradient Mapping:}\\
One source of nonlinearity that commonly arises in optimization algorithms is the gradient $\nabla f(p)$ of a function $f:\mathbb{R}^n\rightarrow\mathbb{R}$. In particular, characterizing the input-output properties of the gradient of a \emph{strongly convex} function with a \emph{Lipschitz continuous} gradient is of interest.
\begin{definition}
	The gradient mapping $\nabla f: \mathbb{R}^n \rightarrow \mathbb{R}^n$ is \textit{Lipschitz continuous} with parameter $L$ if 
	\[
	||\nabla f(p)-\nabla f (p_*)|| \leq L ||p-p_*|| 
	\]
	holds for all $(p,p_*)\in \mathbb{R}^n \times \mathbb{R}^n$, where $L \geq 0$ is a real constant.
\end{definition} 

\begin{definition}
	The function $f: \mathbb{R}^n \rightarrow \mathbb{R}^n$ is said to be \textit{strongly convex} if 
	\[
		(\nabla f(p)-\nabla f(p_*))^T(p-p_*) \geq m ||p-p_*||^2 
	\]
	holds for all $(p,p_*) \in \mathbb{R}^n \times \mathbb{R}^n$, where $m >0$ is a real constant.
\end{definition}

Using these two properties, it is possible to show that $\nabla f$ satisfies the pointwise IQC $(Q_f,p_*, \nabla f(p_*))$ defined by the matrix
\begin{equation} \label{IQCGradient}
	Q_f:=\begin{bmatrix} -2mL & L+m \\ L+m & -2 \end{bmatrix} \otimes I_n.
\end{equation} 
We refer the reader to \cite{Lessard2016} or \cite{Nesterov2013} for a proof of this statement.

\vspace{1ex}
\noindent
\emph{Proximal Mapping:}\\
The second type of nonlinearity that will be used in this paper arises from the proximal mapping of a function.
\begin{definition}
	The proximal mapping $\Pi_{\rho f}: \mathbb{R}^n \rightarrow \mathbb{R}^n$ of the function $f: \mathbb{R}^n \rightarrow \mathbb{R}$ with real parameter $\rho >0$ is defined as
	\begin{equation} \label{eq:proximal-operator}
	\Pi_{\rho f} (p) := \underset{v \in \mathbb{R}^n}{\text{arg min}} \ f(v) + \frac{1}{2 \rho } ||v-p||^2.
	\end{equation}
\end{definition}

The optimality condition of the optimization problem associated with (\ref{eq:proximal-operator}) is:
\begin{equation} \label{proxOPT}
0=\nabla f( \Pi_{\rho f}(p))+\frac{1}{\rho}(\Pi_{\rho f}(p)-p).
\end{equation}
From (\ref{proxOPT}), the proximal mapping can be viewed as the composition of the gradient mapping with an affine operator, followed by an inversion operation:
\[
\Pi_{\rho f}(p) = (I+\rho \nabla f)^{-1}(p).
\]
A known result is that $\Pi_{\rho f}$ satisfies the pointwise IQC $(Q_{\Pi_{\rho f}},p_*,\Pi_{\rho f}(p_*))$ defined by the matrix
\begin{equation*}
	Q_{\Pi_{\rho f}} :=
	\Big(\begin{bmatrix} 0 & \rho^{-1} \\ 1 & -\rho^{-1} \end{bmatrix} \otimes I_n \Big) 
	Q_f 
	\Big( \begin{bmatrix} 0 & 1 \\ \rho^{-1} & -\rho^{-1} \end{bmatrix} \otimes I_n \Big).
\end{equation*}
This result can be derived by using Lemma \ref{AffineIQC} followed by an IQC for inversion operations \cite{Fazlyab2017}.\vspace{.5ex}

\vspace{1ex}
\noindent
\emph{Affine Composition of IQCs:}\\
The following lemma, whose proof can be found in~\cite{Fazlyab2017}, shows how to derive IQCs when a nonlinearity $\phi$ is composed with an affine map.
\begin{lemma}  \label{AffineIQC}
	(IQC for Affine Operations) Consider the nonlinear mapping $\phi$ that satisfies the pointwise IQC defined by $(Q_\phi,p_*, \phi(p_*))$. Define the affine mapping $\psi: \mathbb{R}^n \rightarrow \mathbb{R}^n$ to be $\psi(p):= S_2 p + S_1 \phi (S_0 p)$ where $S_0, S_1, S_2 \in \mathbb{R}^{n \times n}$ and $S_1$ is invertible. Then, $\psi$ satisfies the pointwise IQC defined by $(Q_\psi,p_*, \psi(p_*))$, where
	\begin{equation*}
	Q_\psi := \begin{bmatrix} S_0^T & -(S_1^{-1}S_2)^T \\ 0 & (S_1^{-1})^T \end{bmatrix} 
	Q_\phi
	\begin{bmatrix} S_0 & 0 \\ -S_1^{-1}S_2 & S_1^{-1} \end{bmatrix}.
	\end{equation*}\vspace{0ex}
\end{lemma}

\noindent
\emph{Stability Analysis Using IQCs:}

The following lemma is useful when deriving stability conditions in terms of a LMI. See \cite{Boyd1994} and \cite{Derinkuyu2006} for details.

\begin{lemma} \label{SLemma}
	(Lossless S-Lemma) Let $A^T=A \in \mathbb{R}^{n \times n}$ and $B^T=B \in \mathbb{R}^{n \times n}$. Then, $A \succeq \sigma B$ holds for some $\sigma \geq 0$ if and only if $x^TBx \geq 0 \implies x^TAx \geq 0$ for all $x \in \mathbb{R}^n$. 

We will now show how the input-output properties of an IQC can be used to generate a sufficient stability condition for the feedback interconnection of a LTI system and nonlinearity $\phi$.

\end{lemma}
\begin{proposition}\label{prop:stability}
	Consider a LTI system defined by the matrices $\hat{A} \in \mathbb{R}^{n \times n}$, $\hat{B} \in \mathbb{R}^{n \times m}$, $\hat{C} \in \mathbb{R}^{m \times n}$, and $\hat{D} \in \mathbb{R}^{m \times m}$ with state $\xi \in \mathbb{R}^{n}$, input $ q \in \mathbb{R}^m$, and output $p \in \mathbb{R}^m$:
	\begin{gather*} \label{eq:LTIProp}
	\begin{split}
		\dot{\xi}(t)&=\hat{A}\xi(t)+\hat{B}q(t) \\
		p(t)&=\hat{C}\xi(t)+\hat{D}q(t).
	\end{split}
	\end{gather*}
	Suppose the LTI system has the nonlinearity $\phi: \mathbb{R}^{m} \rightarrow \mathbb{R}^{m}$ as feedback so that $q=\phi(\hat{C}\xi+\hat{D}q)$. Assume $\phi$ satisfies the pointwise IQC ($Q_{\phi}, p_*, \phi(p_*)$) and the feedback interconnection is well-posed.\footnote{The definition of well-posedness is given in Section \ref{sec:problem setup}.} Then, the closed-loop equilibrium point $\xi_* \in \mathbb{R}^{n}$ has global exponential asymptotic stability of at least rate $\alpha$ if the LMI
	\begin{equation} \label{LMI1}
		\begin{bmatrix}
		 \hat{A}^TP + P\hat{A} + \alpha P & P\hat{B} \\ \hat{B}^TP & 0 
		\end{bmatrix} +
		\sigma \begin{bmatrix} \hat{C}^T & 0 \\ \hat{D}^T  & I_{n} \end{bmatrix} Q_{\phi} \begin{bmatrix} \hat{C} & \hat{D} \\ 0 & I_{n} \end{bmatrix} \preceq 0
	\end{equation}
	is feasible for some $\sigma \geq 0$, $\alpha >0$, and $P \succ 0$.
\end{proposition}
\begin{proof}
	Assume that (\ref{LMI1}) is feasible. Let $\delta \xi \!:= \! \xi -  \xi_*$  and $ \delta q \!:= \! q-q_*$, where $q_*$ is the input that achieves equilibrium. Consider the quadratic function $V(\delta \xi)=(\delta \xi)^T P \delta \xi$, where $P \in \mathbb{R}^{n \times n}$, $P \succ 0$. Lyapunov theory states that if $V$ satisfies:
	\begin{itemize}
		\item $V(0)=0$ and $V(\delta \xi)>0$ for all $\delta \xi \in \mathbb{R}^n \setminus \{0\}$,
		\item if $||\delta \xi|| \to \infty$, then $V(\delta \xi) \to \infty$ (radially unbounded),
		\item $\dot{V}(\delta \xi) \leq  -\alpha V(\delta \xi)$ for all $\delta \xi \in \mathbb{R}^n \setminus \{0\}$ and $\alpha >0$,
	\end{itemize}
 	then the equilibrium point has global exponential asymptotic stability of at least rate $\alpha$ \cite{Kahlil2002}.  

	Clearly, $V(0_n)=0_n^T P 0_n=0$. The property $V(\delta \xi)>0$ holds $\forall \delta \xi \neq 0$ because $P \succ 0$. The radial unboundedness property similarly follows from $P \succ 0$. Using the fact that $\hat{A}\xi_*+\hat{B}q_*=0$, the third property can be expressed as
	\begin{gather} \label{stabilityEquality}
	\begin{split}
		 & \dot{V}(\delta \xi) \! + \! \alpha V(\delta \xi) \! = \! 2(\delta \xi)^T P \dot{\xi} \! +  \! \alpha (\delta \xi)^T P \delta \xi  \\
		  &=2(\delta \xi)^T P((\hat{A}\xi \!+ \! \hat{B}q) \! - \! (\hat{A}\xi_* \! + \! \hat{B}q_*)) \! + \! \alpha (\delta \xi)^T P \delta \xi \\
		  & =2(\delta \xi)^T P (\hat{A} \delta \xi \!+ \! \hat{B} \delta q)  \!+ \!  \alpha (\delta \xi)^T P \delta \xi \\
		  & =(\delta \xi)^T P (\hat{A} \delta \xi \! + \! \hat{B}\delta q) \! + \!(\hat{A}\delta \xi \!+ \! \hat{B} \delta q)^T P \delta \xi \! + \! \alpha (\delta \xi)^T P \delta \xi \\  
		  & = \begin{bmatrix}
		  \delta \xi\\ \delta q
		  \end{bmatrix} ^T \begin{bmatrix}
		  	\hat{A}^T P \! + \! P\hat{A} \! + \! \alpha P & P\hat{B} \\\hat{B}^T P & 0
		  \end{bmatrix}
		  \begin{bmatrix}
		  \delta \xi \\ \delta q
		  \end{bmatrix} \leq 0.		 
	\end{split}
	\end{gather}
	Finally, since the pointwise IQC $(Q_{\phi},p_*, \phi(p_*)$ is satisfied,
	\begin{equation} \label{pointwiseIQC}
		\begin{bmatrix} \delta \xi \\ \delta q \end{bmatrix}^T \begin{bmatrix} \hat{C}^T & 0 \\ \hat{D}^T & I_{n} \end{bmatrix} Q_\phi \begin{bmatrix} \hat{C} & \hat{D} \\ 0 & I_{n} \end{bmatrix} \begin{bmatrix} \delta \xi \\ \delta q \end{bmatrix} \geq 0.
	\end{equation}
	Since (\ref{LMI1}) is feasible and (\ref{pointwiseIQC}) holds, it directly follows from Lemma \ref{SLemma} that (\ref{stabilityEquality}) holds. Hence, the equilibrium $\xi_*$ has global exponential asymptotic stability of at least rate $\alpha$.
\end{proof}

\tomove{
\enrique{Text below contains description of the  results. Will move to another section.}

\subsection{IQC for the Gradient Mapping}
One design for the optimization mapping would be the negative gradient mapping $\Psi_1= - \nabla f$ because the output is zero if and only if the input is an optimal solution of (\ref{OPT}). Note that the negative gradient is considered rather than the gradient for stability purposes, as will be discussed later. In addition to $f$ being a strongly convex function, it is assumed that $\nabla f$ is Lipschitz continuous. Definitions of these assumptions are now stated.

\begin{lemma} \label{IQCNegGradient}
	Assume the pointwise IQC $(Q_f, p_*, \nabla f(p_*))$ is satisfied and let $\Psi_1 = - \nabla f$. Then, the pointwise IQC $(Q_{\Psi_1}, p_*, \Psi_1(p_*))$ defined by the matrix
	\begin{equation*}
		Q_{\Psi_1} = \Big ( \begin{bmatrix} 1 & 0 \\ 0 & -1 \end{bmatrix} \otimes I_n \Big ) Q_f \Big (  \begin{bmatrix} 1 & 0 \\ 0 & -1 \end{bmatrix} \otimes I_n \Big ).
 	\end{equation*} 
 	is satisfied.
\end{lemma}
\begin{proof}
	The IQC immediately follows from Lemma \ref{AffineIQC} with $\phi= \nabla f$, $S_0=I_n$, $S_1=-I_n$, and $S_2=0_{n \times n}$. 
\end{proof}

Another possible design choice for the optimization mapping would be to incorporate the proximal mapping by choosing $\Psi_2 = \Pi_{\rho f} - I_n$. The following proposition shows why this is a satisfactory choice. 

\begin{proposition}
The output of the mapping $\Psi_2$ is zero if and only if the input is an optimal solution of (\ref{OPT}).
\end{proposition}
\begin{proof}
	($\implies$): Assume that $\Psi_2(p_*)=0$. Then from the definition of $\Psi_2$, $\Pi_{\rho f} (p_*) =  p_*$. In terms of the proximal mapping definition, $q=p_*$. It follows from the optimality condition of the proximal mapping optimization problem that
	\begin{equation*}
		0 = \nabla f(q) + \frac{1}{\rho}(q-p_*) \implies 0=\nabla f(p_*) \implies p_* \in \mathcal{X}^*.
	\end{equation*}
	($\impliedby$): Assume that $p_* \in \mathcal{X}^*$. It follows from 
	\begin{equation*}
		\underset{q \in \mathbb{R}^n}{\text{arg min}} \ f(q)=p_* \text{ and } \underset{q \in \mathbb{R}^n}{\text{arg min}} \ \frac{1}{2 \rho} ||q-p_*||^2 = p_*
	\end{equation*}
	that $\Pi_{\rho f}(p_*)=p_*$. Thus, $\Psi_2(p_*)=0$.
\end{proof}

\begin{lemma} \label{IQCOptimization}
	Assume the pointwise IQC $(Q_f, p_*,\nabla f(p_*))$ is satisfied and let $\Psi_2=\Pi_{\rho f}-I_n$. Then, the pointwise IQC $(Q_{\Psi_2},p_*,\Psi_2(p_*))$ defined by the matrix 
	\begin{equation*}
		Q_{\Psi_2}=
		\Big ( \begin{bmatrix} 1 & 1 \\ 0 & 1 \end{bmatrix} \otimes I_n \Big )
		Q_{\Pi_{\rho f}}
		\Big ( \begin{bmatrix} 1 & 0 \\ 1 & 1 \end{bmatrix} \otimes I_n \Big ).
	\end{equation*}
	is satisfied.
\end{lemma}
\begin{proof}
	The IQC immediately follows from Lemma \ref{AffineIQC} with $\phi=\Pi_{\rho f}$, $S_0=I_n$, $S_1=I_n$, and $S_2=-I_n$.
\end{proof}
}

\section{Problem Setup}\label{sec:problem setup}

\begin{figure}[h]
	\centering
	\includegraphics[width=0.65\linewidth]{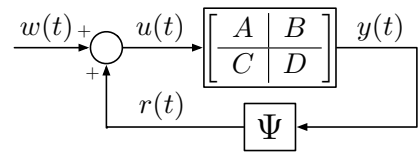}
	\caption{LTI system interconnected with a nonlinear mapping and constant disturbance signal.} 
	\label{fig:system}
\end{figure}

The problem setup is illustrated in Fig. \ref{fig:system}, where we consider a LTI system represented by a state-space model where $x \in \mathbb{R}^n$ is the state, $u \in \mathbb{R}^m$ is the input, and $y \in \mathbb{R}^p$ is the output:
\begin{gather}
\begin{split}
& \dot{x}(t)=Ax(t)+B u(t)  \label{eq:LTI}  \\
& y(t)=Cx(t)+D u(t).
\end{split}
\end{gather}
The input $u(t)$ is the sum of a control signal $r(t)\in \mathbb{R}^m$ and an unknown constant disturbance $w(t)=w\in\mathbb{R}^m$, i.e. $u(t)=w(t)+r(t)=w+r(t)$.
Finally, the feedback operator $\Psi(\cdot)$ denotes the (possibly nonlinear) feedback control to be designed.

Our goal is to, given the measurement $y$, design a control input $r=\Psi(y)$ that drives the system \eqref{eq:LTI} to a steady-state $x_*$ that is an optimal solution of a predefined optimization problem
\begin{gather} \label{OPT}
	\min_{x \in \mathbb{R}^n} f(x), 
\end{gather}
where $f: \mathbb{R}^n  \rightarrow \mathbb{R}$ is a given cost function.

Therefore, given the measurement $y(t)$, the feedback $\Psi(\cdot)$ must produce a control $r = \Psi(y)$ such that $x(t)\rightarrow \mathcal X^*$, where $\mathcal X^*$ is the set of optimal solutions to (\ref{OPT}), i.e.,
\begin{equation*} 
	\mathcal{X}^* =\{ x \in \mathbb{R}^n : \nabla f (x)=0 \}. 
\end{equation*}
Throughout this paper we make the following assumption.
\begin{assumption}
The cost function $f(x)$ of the optimization problem (\ref{OPT}) is continuously differentiable, strongly convex, and has a Lipschitz continuous gradient. This implies that the set $\mathcal X^*$ is a singleton.\footnote{Relaxing this assumption is desired and is a subject of future work.}
\end{assumption}


Finally, we provide a concrete model for $\Psi(\cdot)$. As the optimality conditions for optimization problem \eqref{OPT} are in general nonlinear, the feedback controllers to be designed will be necessarily of the same type. Thus, we consider the nonlinear feedback $\Psi$ using the nonlinear dynamics
\begin{gather} \label{NonlinearDynamics}
\Psi: \quad \begin{split}
&\dot{\eta}(t)=F(\eta(t),y(t)) \\
& r(t)=H(\eta(t),y(t)),
\end{split}
\end{gather}
where $\eta \in \mathbb{R}^d$ is the state of the feedback dynamics, $r \in \mathbb{R}^m$ is the output of the feedback dynamics, and the mappings $F: \mathbb{R}^d \times \mathbb{R}^p \rightarrow \mathbb{R}^d$ and $H: \mathbb{R}^d\times\mathbb{R}^p \rightarrow \mathbb{R}^m$ are possibly nonlinear.

\begin{remark}[Well-Posedness]
From the feedthrough terms present in \eqref{eq:LTI} and \eqref{NonlinearDynamics}, it is possible a priori that the feedback interconnection is not well-posed.\footnote{The feedback interconnection of \eqref{eq:LTI} and \eqref{NonlinearDynamics} is well-posed if $u(t)$ and $y(t)$ are uniquely defined for every choice of states $x(t)$ and $\eta(t)$.}
A sufficient condition that prevents this problem is by enforcing that whenever $D\not=0$, the map $H$ depends only on $\eta$, i.e., $r(t)= H(\eta(t))$. We will further discuss this condition in Section \ref{sec:control design}.
\end{remark}

\subsection{Design Challenges}
There are several challenges associated to designing \eqref{NonlinearDynamics} such that in steady state $x_*\in\mathcal X^*$.
\begin{itemize}
\item \emph{Lack of direct access to $x(t)$:} The system output matrix $C$ is not necessarily invertible. Thus, recovering $x(t)$ from $y(t)$ is not straightforward. 
\item \emph{Finding the solution $x_*\in\mathcal X^*$:} Finding the optimal solution to the optimization problem is usually challenging or the cost function may change, giving not enough time to recompute $x_*$.
\item \emph{Driving $x(t)$ to  $x_*\in\mathcal X^*$:} Even if one has access to the optimal solution $x_*$, one then needs to design the right $r(t)$ that ensures that $x(t)$ converges to it.
\end{itemize}

Interestingly, some of these challenges can be easily handled using tools from control theory, such as recovering $x(t)$ from $y(t)$ or driving $x(t)$ to $x_*$. On the other hand, finding an optimal solution $x_*$ is the major goal within optimization theory. Therefore, when the timescale of the control and optimization tasks do not intersect, our problem can be easily solved using standard tools from control and optimization. However, when the timescale separation is no longer present, the problem becomes more challenging as there is no standard tools to address it. In particular, it is usually hard to assess the stability of such an interconnected system. This problem is systematically addressed in the next section.

\section{Optimization-Based Control Design}\label{sec:control design}

In this section we describe the proposed  optimization-based controllers, that combine tools from control and optimization, and leverage the IQC framework described in the preliminaries. The crux of our solution is a modularized architecture that breaks down the feedback dynamics (\ref{NonlinearDynamics}) into three serial components that systematically addresses the challenges described in the previous section and allow a straightforward application of Proposition \ref{prop:stability} to certify global exponential asymptotic stability.
\begin{figure}[!htb]
	\centering
	\includegraphics[width=.65\linewidth]{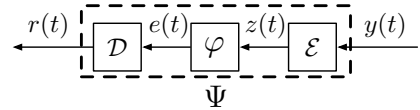}
	\caption{Optimization-based control feedback breakdown.} 
	\label{fig:system2}
\end{figure} 

The proposed architecture is described in Fig. \ref{fig:system2}.
The first component $\mathcal E: y(t)\mapsto  z(t)$ is a state \emph{estimator} that takes the output of the LTI system and produces a state estimate $ z(t)$. The second component $\varphi: z(t)\mapsto e(t)$, referred to as the \textit{optimizer}, takes the state estimate and produces a measurement of the optimality error or direction of desired drift $e(t)$, which is required to be zero if and only if the input is in the set $\mathcal{X}^*$. The optimizer can be thought of as the part of optimization algorithm that dictates the direction of the next step. The third component $\mathcal D: e(t)\mapsto r(t)$, the \emph{driver}, takes the optimality error and produces the input to the LTI system that ensures that the equilibrium satisfies $e_*= 0$.

\begin{remark}
One of the advantages of the proposed architecture is the role independence of each component. This allows for a subsystem to be skipped if the functionality is not required. 
For example, in cases where $y(t)=x(t)$ or the optimization problem uniquely depends on $y(t)$, then the estimator block can be avoided.
\end{remark}

In the remainder of this section, we describe the design requirements of each proposed component/subsystem and give some examples on how to implement them.

\subsection{Design of the Estimator $\mathcal E$}

The estimator component $\mathcal E: y(t)\mapsto z(t)$ is perhaps the simplest to design. Its goal to build an estimate of the state, $z(t)$, from $y(t)$. The design requirement of $\mathcal E$ is:
\begin{itemize}
	\item \textbf{A.1:} If the system is in equilibrium, $z(t)=z_* =x_*$.
\end{itemize}

Therefore, an obvious choice for $\mathcal E$ is an observer/state estimator. 
The dynamics of $\mathcal E$ are therefore given by
\begin{gather*} 
\mathcal{E}: \quad \begin{split}
&\dot{\hat{x}} = (A-LC)\hat{x}+(B-LD)u + Ly \\
& z=\hat{x},
\end{split}
\end{gather*}
where $L \in \mathbb{R}^{n \times p}$ is a constant matrix to be designed.

A standard argument for observers shows that the evolution of the error $\delta x(t) := x(t)-z(t)$ is given by
\[
\dot{\delta x}(t) = (A-LC)\delta x(t).
\]
Moreover, if \eqref{eq:LTI} is observable, $L$ can be chosen to satisfy
\begin{equation} \label{observable}
	\text{rank}(A-LC)=n.
\end{equation}

\subsection{Design of the Optimizer $\varphi$}

The optimizer $\varphi$ has two design requirements:
\begin{itemize}
	\item \textbf{B.1:} The optimizer must take the estimated state $z(t)$ as an input and then produce a measure of optimality error or direction of drift $e(t)$ such that $e(t)=0$ if and only if $z(t)=x_*\in\mathcal{X}^*$.
	\item \textbf{B.2:} The input-output characteristics of $\varphi$ must be captured by an IQC $(Q_\varphi, z_*,\varphi(z_*))$.
\end{itemize}

For the purpose of this paper, we consider two possible solutions.\vspace{1ex}

\noindent
\emph{$\varphi_1$: Gradient Descent.}~
The first solution considered is the standard gradient descent mapping, i.e., 
\begin{equation} \label{eq:gradient-descent}
	\varphi_1 := -\nabla f.
\end{equation}
It is straightforward to verify that $e(t)=-\nabla f(z(t))=0$ if and only if $z(t)=x_*\in\mathcal X^*$. 
The following lemma explicitly computes the IQC for $\varphi_1$.

\begin{lemma} \label{IQCNegGradient}
	Assume the pointwise IQC $(Q_f, z_*, \nabla f(z_*))$ is satisfied. Then, the pointwise IQC $(Q_{\varphi_1}, z_*, \varphi_1(z_*))$ defined by the matrix
	\begin{equation*}
		Q_{\varphi_1} := \Big ( \begin{bmatrix} 1 & 0 \\ 0 & -1 \end{bmatrix} \otimes I_n \Big ) Q_f \Big (  \begin{bmatrix} 1 & 0 \\ 0 & -1 \end{bmatrix} \otimes I_n \Big )
 	\end{equation*} 
 	is satisfied.
\end{lemma}
\begin{proof}
	The IQC immediately follows from Lemma \ref{AffineIQC} with $\phi= \nabla f$, $S_0=I_n$, $S_1=-I_n$, and $S_2=0_{n \times n}$. 
\end{proof}\vspace{1ex}

\noindent
\emph{$\varphi_2$: Proximal Tracking.}~
Our second option for the optimizer block is inspired by the proximal mapping \eqref{eq:proximal-operator}.
It essentially computes the error between the input $z(t)$ and the solution given by the proximal operator $\Pi_{\rho f}(z(t))$, that is, 
\begin{equation}\label{eq:proximal-tracking}
 \varphi_2 :=\Pi_{\rho f}-I_n.
\end{equation}
The following proposition shows that \eqref{eq:proximal-tracking} satisfies the first design requirement.
\begin{proposition}
The mapping $\varphi_2$ satisfies the property that $e(t)=\Pi_{\rho f}(z(t))-z(t)=0$ if and only if $z(t)=z_* \in \mathcal{X}^*$.
\end{proposition}
\begin{proof}
	Let $z(t)=z_*$ and assume that $\varphi_2(z_*)=0$. Then from (\ref{eq:proximal-tracking}), $\Pi_{\rho f} (z_*) =  z_*$. It follows from (\ref{proxOPT}) that
	\begin{equation*}
		0 = \nabla f(\Pi_{\rho f}(z_*)) + \frac{1}{\rho}(\Pi_{\rho f}(z_*)-z_*) \iff 0=\nabla f(z_*).
	\end{equation*}
	By the definition of $\mathcal{X}^*$, $z_* \in \mathcal{X}^*$.
	
	Conversely, assume that $z_* \in \mathcal{X}^*$. Since
	\begin{equation*}
		\underset{v \in \mathbb{R}^n}{\text{arg min}} \ f(v)=z_* \text{ and } \underset{v \in \mathbb{R}^n}{\text{arg min}} \ \frac{1}{2 \rho} ||v-z_*||^2 = z_*, 
	\end{equation*}
	\[
	\text{we have } \underset{v \in \mathbb{R}^n}{\text{arg min}} \ f(v) + \frac{1}{2 \rho} ||v-z_*||^2 = z_*.
	\]	
	This is equivalent to $\Pi_{\rho f}(z_*)=z_*$. Thus, $\varphi_2(z_*)=0$.
\end{proof}

Finally, the next lemma computes the IQC that characterizes $\varphi_2$.
\begin{lemma} \label{IQCOptimization}
	Assume the pointwise IQC $(Q_f, z_*,\nabla f(z_*))$ is satisfied. Then, the pointwise IQC $(Q_{\varphi_2},z_*,\varphi_2(z_*))$ defined by the matrix 
	\begin{equation*}
		Q_{\varphi_2} :=
		\Big ( \begin{bmatrix} 1 & 1 \\ 0 & 1 \end{bmatrix} \otimes I_n \Big )
		Q_{\Pi_{\rho f}}
		\Big ( \begin{bmatrix} 1 & 0 \\ 1 & 1 \end{bmatrix} \otimes I_n \Big ).
	\end{equation*}
	is satisfied.
\end{lemma}
\begin{proof}
	The IQC immediately follows from Lemma \ref{AffineIQC} with $\phi=\Pi_{\rho f}$, $S_0=I_n$, $S_1=I_n$, and $S_2=-I_n$.
\end{proof}

\subsection{Design of the Driver $\mathcal D$}
The last component of the proposed solution is in charge of generating the control signal $r(t)$ that drives the system towards the optimal solution of \eqref{OPT}. The design requirement for $\mathcal{D}$ is:
\begin{itemize}
	\item \textbf{C.1:} If the system is in equilibrium, then $e(t)=e_*=0$.
\end{itemize}

Thus, one possible choice would be to use a Proportional-Integral (PI) controller defined by the dynamics
\begin{gather*} 
\mathcal{D}: \quad \begin{split}
& \dot{e}_I= e, \ r= K_I e_I+K_Pe
\end{split}
\end{gather*}
where $K_I, K_P \in \mathbb{R}^{m \times n}$ are constant matrices to be designed. 

It is straightforward to show that $\dot e_I=0$ if and only if $e(t)=0$, which satisfies our design requirements. In fact, this also shows that we only need an integrator to satisfy the design requirement. However, a PI  controller provides better dynamic properties than a pure integrator and therefore we choose to add the proportional term.

\subsection{Integrated System}

\begin{figure}[!htb]
	\centering
	\includegraphics[width=1\linewidth]{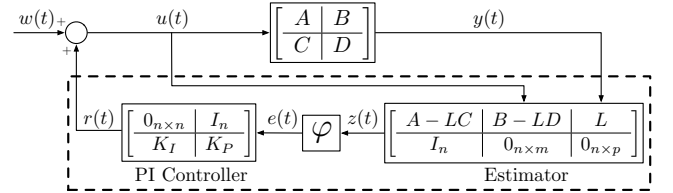}
	\caption{LTI system in feedback with a state estimator $\mathcal{E}$, optimizer $\varphi$, and PI controller as the driver $\mathcal{D}$.} 
	\label{fig:case2system}
\end{figure}

This resulting interconnected system is shown in Fig. \ref{fig:case2system}.

In terms of the dynamics defined by (\ref{NonlinearDynamics}), the feedback design is given by
\begin{align}
& \dot{\eta} \!  = \!\begin{bmatrix} A\! -LC \! + \!(B-LD) K_P \varphi & (B \! -LD) K_I \\ \varphi & 0 \nonumber\end{bmatrix}  \eta \! + \! \begin{bmatrix} L \\ 0 \end{bmatrix}y   \\
 &r  =\begin{bmatrix} K_P\varphi & K_I \end{bmatrix}\eta, \text{ where } \eta := \begin{bmatrix} \hat{x}^T & e_I^T \end{bmatrix}^T \in \mathbb{R}^{2n}. \label{NonlinearDynamicsChoice}
\end{align}
Here we have used the notation $\varphi \hat{x} := \varphi(\hat{x})$.

\begin{remark}
Whenever the estimator subsystem $\mathcal{E}$ is included in the interconnection, the feedback interconnection will be well-posed because $r(t)$ only depends on $\eta(t)$. However, well-posedness is not guaranteed when $D \neq 0$, $K_P \neq 0$, and there is no estimator subsystem because
\begin{equation*}
	r(t)=K_P \varphi (Cx(t)+D(w+r(t))) +K_I e_I(t)
\end{equation*}
 depends on itself. One simple solution to overcome this issue is add a module $\mathcal{E}$ such that $z(t)=y(t)-Du(t)$.
\end{remark}

The next section shows that indeed the integrated system is able to guarantee steady-state optimality under mild conditions and illustrates how the IQC framework can be leveraged to guarantee global exponential asymptotic stability.

\section{Optimality and Convergence}\label{sec:analysis}

\subsection{Optimality Analysis}
The optimality analysis requires the following assumption.
\begin{assumption}\label{as:B-inv}
	The system is steady-state controllable. That is, given any steady-state $x_*$, there exists an input $u_*$ such that $Ax_*+Bu_*=0$.
\end{assumption}

Assumption \ref{as:B-inv} is in some sense \emph{necessary} to ensure that the system can achieve an arbitrary steady-state. Although this assumption is stronger than the standard controllability assumption, we point out that while controllability is sufficient to drive $x(t)$ towards any state $x_*$ in finite time, it does not requires that $x(t)$ remains equal to $x_*$.

\begin{theorem} \label{Optimality}
	Consider the interconnection of the LTI system (\ref{eq:LTI}) and nonlinear feedback (\ref{NonlinearDynamicsChoice}), where design requirements A.1, B.1, B.2, and C.1 are satisfied. Suppose (\ref{eq:LTI}) is a minimal realization and Assumption 2 is satisfied. Then, $(x_*,\eta_*)$ is an equilibrium point of the interconnected system for some point $\eta_*$ and matrix $L$ that satisfies (\ref{observable}) if and only if $x_* \in \mathcal{X}^*$.
\end{theorem}
\begin{proof}
	Assume $(x_*,\eta_*)$ is an equilibrium point of the interconnected system, where $\eta_* = \begin{bmatrix} z^T_* & e_{I*}^T \end{bmatrix}^T$. The LTI system must then be in equilibrium, meaning that 
	\begin{equation} \label{LTIequil}
		0=Ax_*+Bu_*, \ y_*=Cx_*+Du_*, \ u_*= r_* +w.
	\end{equation}
	Additionally, the dynamics of $\Psi$ must be in equilibrium, meaning that components $\mathcal{E}$ and $\mathcal{D}$ are in equilibrium. It follows from $\mathcal{D}$ being in equilibrium and B.1 that
	\begin{equation} \label{eq:estimatorOptimal}
		\dot{e}_I=e=\varphi(z_*)=0 \iff z_* \in \mathcal{X}^*.
	\end{equation}
	It follows from $\mathcal{E}$ being in equilibrium that
	\begin{gather*}
	\begin{split}
		\dot{\hat{x}}& =(A-LC)z_* +(B-LD)u_* +L(Cx_*+Du_*) \\
		& = (A-LC)z_* +Bu_* +LCx_*=0.
	\end{split}
	\end{gather*}
	Adding $Ax_*$ to both sides and using (\ref{LTIequil}) results in
	\begin{equation*}
		(A-LC)(z_*-x_*)=0.
	\end{equation*}
	Since the LTI system is observable, we can choose $L$ such that $A-LC$ is Hurwitz and therefore $x_* = z_*$. From (\ref{eq:estimatorOptimal}), $ x_* \in \mathcal{X}^*$. 

	Conversely, assume that $x_* \in \mathcal{X}^*$. Consider $e_*=0$ and an $e_{I*}$ such that $Ax_*+B(K_Ie_{I*}+w)=0$, which exists because of Assumption 2. Then, $\dot{x}=0$ and $\dot{e}_I=0$ directly follow. Next, consider $z_*=x_*$. It is then straightforward to show that $\dot{\hat{x}}=0$. The optimality property of $\varphi$ gives that $\varphi(z_*)=e_*=0$, which is consistent with the previous definition of $e_*$. Therefore, $(x_*,\eta_*)$ is an equilibrium point.
\end{proof}
\subsection{Stability Analysis} 
This section will derive a sufficient condition for the global exponential asymptotic stability of the equilibrium point considered in the optimality analysis. For this analysis, it is useful to group the linear dynamics of $\mathcal{E}$ and $\mathcal{D}$ into the LTI system to essentially create a larger dimension LTI system. The resulting system is 
\begin{align*}
	\dot{\xi} &= \underbrace{\begin{bmatrix} A & 0 & B K_I \\ LC & A-LC & B K_I \\ 0 & 0 & 0  \end{bmatrix}}_{\hat{A}} \xi
	+ \underbrace{\begin{bmatrix} B K_P \\ B K_P  \\ I_n  \end{bmatrix}}_{\hat{B}_e} 
	e
	+ \underbrace{ \begin{bmatrix} B \\ B \\ 0 \end{bmatrix}}_{\hat{B}_w} w \\
	z&= \underbrace{\begin{bmatrix} 0 & I_n & 0 \end{bmatrix}}_{\hat{C}} \xi, \text{where $\xi := \begin{bmatrix} x^T & \hat{x}^T & e_I^T \end{bmatrix}^T$.}
\end{align*}

\begin{theorem}\label{th:global-stability} 
	Consider the interconnection of the LTI system (\ref{eq:LTI}) and nonlinear feedback (\ref{NonlinearDynamicsChoice}) with the equilibrium point considered in Theorem \ref{Optimality}. Assume the pointwise IQC $(Q,z_*, \varphi(z_*))$ is satisfied and the assumptions of Theorem \ref{Optimality} hold. Then, the equilibrium point $(x_*,\eta_*)$ has global exponential asymptotic stability of at least rate $\alpha$ if the LMI 
	\begin{equation} \label{LMI}
	\begin{bmatrix} \hat{A}^TP+P\hat{A}+\alpha P & P\hat{B}_e \\ \hat{B}_e^TP & 0 \end{bmatrix}
	+ \sigma \begin{bmatrix} \hat{C}^T & 0 \\ 0 & I_{n} \end{bmatrix} Q \begin{bmatrix} \hat{C} & 0 \\ 0 & I_{n} \end{bmatrix} \preceq 0
	\end{equation}
	is feasible for some $\sigma \geq 0, \alpha > 0,$ and $ P \succ 0$.
\end{theorem}
\begin{proof}
	This is essentially an application of Proposition \ref{prop:stability}, where the Lyapunov function is given by
	\begin{equation*}
	V(\delta \xi)=(\delta \xi)^TP \delta \xi > 0, \ P \in \mathbb{R}^{3n \times 3n}, \ P \succ 0.
	\end{equation*}
	Using the fact that $\hat{A}z_*+\hat{B}_e e_* + \hat{B}_w w =0$,
	\begin{align*}
	\begin{split}
	\dot{V}(\delta \xi) &+\alpha V(\delta \xi)\\ 
	& = 2(\delta \xi)^TP[\hat{A}\xi+\hat{B}_e e +\hat{B}_w w] + \alpha(\delta \xi)^T P \delta \xi \\
	& = 2(\delta \xi)^TP[\hat{A}\delta \xi+\hat{B}_e\delta e] + \alpha (\delta \xi)^T P \delta \xi \\
	& = \begin{bmatrix} \delta \xi \\ \delta e \end{bmatrix}^T 
	\begin{bmatrix} \hat{A}^T P + P\hat{A}+\alpha P & P \hat{B}_e \\
	\hat{B}_e^T P & 0 \end{bmatrix} \begin{bmatrix} \delta \xi \\ \delta e \end{bmatrix} \leq 0.
	\end{split}
	\end{align*}
	The rest of the proof follows from Proposition \ref{prop:stability}.
\end{proof}

\subsection{Convergence Rate}
Finally, we show how the LMI condition derived in Theorem \ref{th:global-stability} can be leveraged to compute the maximum convergence rate that the system can achieve. Our goal here is to solve the optimization problem:
\begin{equation} \label{alphaOpt}
	\underset{\sigma \geq 0, \alpha >0, P \succ 0 }{\text{maximize}} \ \alpha  \text{ subject to } (\ref{LMI}).
\end{equation}

The main challenge is that because $\alpha$ multiplies $P$ in \eqref{LMI}, the optimization problem is non-convex. However, for a fixed $\alpha$, finding whether \eqref{LMI} is feasible can be done efficiently. Therefore, it is possible to implement a line search in $\alpha$ that finds the maximum value $\alpha_\text{max}$ that satisfies \eqref{LMI}.

\section{Numerical Examples}

\subsection{Scalar System}
Consider the scalar system  
\begin{equation*}
\dot{x}=-5x+u, \ y=x,
\end{equation*}
with $x,u,y \in \mathbb{R}$. Since the estimator module $\mathcal{E}$ is not needed, the dimension of the LMI can be reduced. The feedback interconnection is still well-posed because $D=0$. Let the cost function be of the form
\begin{equation*}
	f(x)=\frac{1}{2}qx^2+cx+v,
\end{equation*}
where $q,c,v \in \mathbb{R}$ are constants. For this case, the Lipschitz constant and strong convexity constant are $m=L=q$. Let the control parameters be given by $k_i=1$ and $k_p=1$.

\begin{figure}[!ht]
	\centering
	\includegraphics[width= \linewidth]{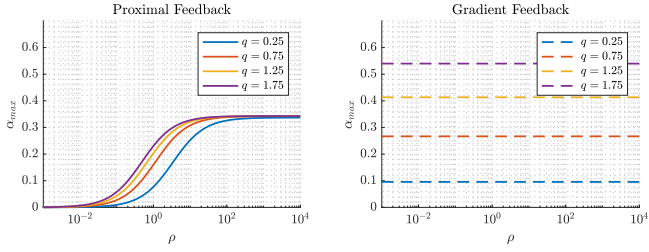}
	\caption{Maximum feasible $\alpha$ versus $\rho$ for several different scalar cost functions when using (a): $\varphi_2$ (proximal optimizer) and (b): $\varphi_1$ (gradient optimizer).}
	\label{fig:example1a}
\end{figure}

Fig. \ref{fig:example1a} shows the solution of (\ref{alphaOpt}) as a function of $\rho$. Several curves representing different steepness levels of $f$ are plotted. The plots demonstrate that larger values of $\rho$ lead to a larger $\alpha_\text{max}$, with no marginal improvement after a certain point. Additionally, it is interesting to note that there are cases when the $\varphi_2$ optimizer achieved a larger $\alpha_\text{max}$ than the $\varphi_1$ optimizer.

\begin{figure}[!ht]
	\centering
	\includegraphics[width= \linewidth]{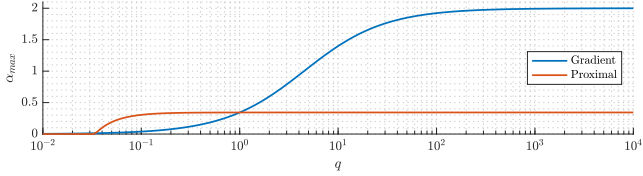}
	\caption{Maximum feasible $\alpha$ versus $q$ when using $\varphi_1$ and $\varphi_2$ optimizers.}
	\label{fig:example1b}
\end{figure}

After choosing a sufficiently large $\rho$, the solution of (\ref{alphaOpt}) was plotted as a function of $q$ as shown in Fig. \ref{fig:example1b}. As expected, larger values of $q$, corresponding to steeper quadratic functions, resulted in a larger $\alpha_\text{max}$. There was also no marginal improvement past a certain threshold of $q$. This threshold was a very large $q$ for $\varphi_1$ and a very small $q$ for $\varphi_2$. For $q<1$, $\varphi_2$ achieved a larger $\alpha_\text{max}$ and for $q>1$, $\varphi_1$ achieved a larger $\alpha_\text{max}$.

\begin{figure}[!ht]
	\centering
	\includegraphics[width= \linewidth]{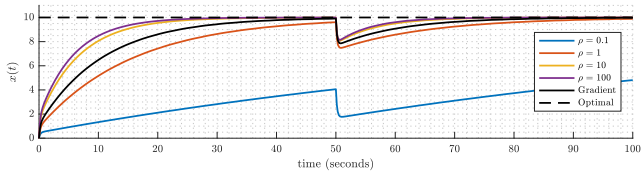}
	\caption{System state versus time when using $\varphi_1$ and $\varphi_2$ optimizers. Several choices of $\rho$ are shown for the $\varphi_2$ optimizer.}
	\label{fig:example1c}
\end{figure}

The system's state as a function of time, when the cost function was $f(x)=(x-10)^2$, is given in Fig. \ref{fig:example1c}. Several trajectories, corresponding to different values of $\rho$, were plotted. The disturbance signal was initially set as $w=2$ and at $t=50$s was changed to $w=-10$. For all cases, the state was able to recover from the change in disturbance and continue tracking the optimal solution. The trajectories illustrate that the performance of the $\varphi_2$ optimizer is very much related to the choice of $\rho$ and if chosen correctly can outperform the $\varphi_1$ optimizer.  

\subsection{MIMO System with State Estimator}
Consider the MIMO system defined by
\begin{equation*}
	A=\begin{bmatrix} 0 & 1 \\ -10 & -5 \end{bmatrix}, \ B= \begin{bmatrix} 1  & 4 \\ 1 & 0 \end{bmatrix}, \ C = \begin{bmatrix} 1 & 0 \end{bmatrix}, \ D=0_{1 \times 2},
\end{equation*}
with $x,u \in \mathbb{R}^2$ and $y \in \mathbb{R}$. Since the output of the LTI system only has information about the first state, an estimator module is obviously needed. Let the cost function be of the form
\begin{equation*} 
	f(x)=\frac{1}{2}x^TQx+c^Tx,
\end{equation*}
where $Q \in \mathbb{R}^{2 \times 2}$, $Q \succ 0 $, and $c \in \mathbb{R}^2$. For this case, the Lipschitz constant $L$ is the larger eigenvalue of $Q$ and the strong convexity constant $m$ is the smaller eigenvalue of $Q$. Let the feedback parameters be given by
\begin{equation*}
K_I=K_P= \begin{bmatrix} 0 & 1 \\ 1/4 & -1/4 \end{bmatrix} \text{ and } 
L= \begin{bmatrix} 1 & 1 \end{bmatrix}.
\end{equation*}

\begin{figure}[!ht]
	\centering
	\includegraphics[width= \linewidth]{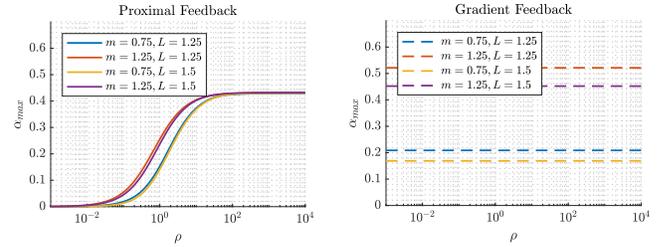}
	\caption{Maximum feasible $\alpha$ versus $\rho$ for several different multivariable cost functions when using (a): $\varphi_2$ (proximal optimizer) and (b): $\varphi_1$ (gradient optimizer).}
	\label{fig:example2a}
\end{figure}

Fig. $\ref{fig:example2a}$ shows the solution of (\ref{alphaOpt}) as a function of $\rho$. Several curves corresponding to different eigenvalue choices of $Q$ are plotted. With a sufficiently large $\rho$, $\varphi_2$ was able to achieve a larger $\alpha_\text{max}$ than $\varphi_1$ when $m=0.75$, but was not able to when $m=1.25$. For both optimizer types, $L=1.25$ resulted in a larger $\alpha_\text{max}$ than $L=1.5$.

\begin{figure}[!ht]
	\centering
	\includegraphics[width= \linewidth]{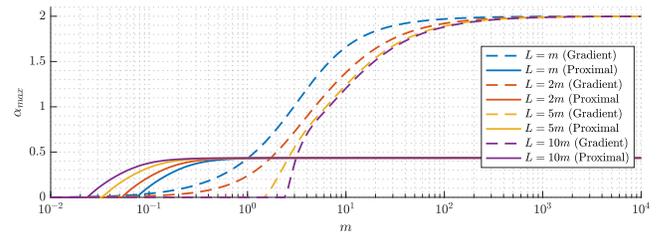}
	\caption{Maximum feasible $\alpha$ versus $m$ for several different choices of $L$ when using $\varphi_1$ and $\varphi_2$ optimizers.}
	\label{fig:example2b}
\end{figure}

After choosing a sufficiently large $\rho$, the solution of (\ref{alphaOpt}) was plotted as a function of $q$ as shown in Fig. \ref{fig:example2b}. The Lipschitz constant was chosen as different multiples of $m$. The $\varphi_1$ optimizer resulted in a larger $\alpha_\text{max}$ when $L$ was chosen closer to $m$. Conversely, the $\varphi_2$ optimizer resulted in a larger $\alpha_\text{max}$ when the multiple was chosen farther from $m$.

\begin{figure}[!ht]
	\centering
	\includegraphics[width= \linewidth]{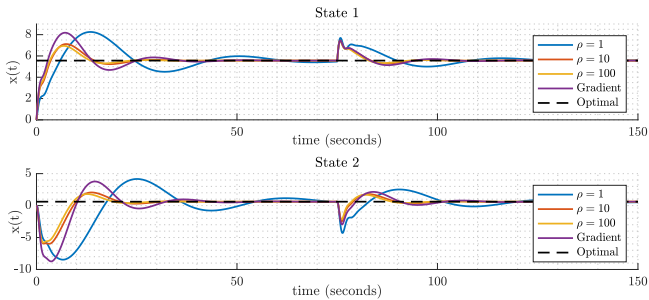}
	\caption{System (a) state 1 and (b) state 2 versus time when using $\varphi_1$ and $\varphi_2$ optimizers. Several choices of $\rho$ are shown for the $\varphi_2$ optimizer.}
	\label{fig:example2c}
\end{figure}

The system's state as a function of time, when the cost function was defined by
\begin{equation*}
	Q= \begin{bmatrix} 1 & 1/6  \\  1/6 & 2/3 \end{bmatrix}, \ c^T= \begin{bmatrix} -17/3 & -4/3 \end{bmatrix},
\end{equation*}
is given in Fig. \ref{fig:example2c}. In this case, $m\approx 0.5976$, $L \approx 1.0690$, and the optimal solution is $x_* \approx \begin{bmatrix} 5.5652 & 0.6087 \end{bmatrix}^T$. The disturbance was initially set to zero, but at $t=75$s was changed to $w= \begin{bmatrix} 1 & 1 \end{bmatrix}^T$. Similar to the scalar case, the trajectories were able to recover from the change in disturbance and there were cases when $\varphi_2$ outperformed $\varphi_1$. It is particularly interesting that state 2 has the ability to reach the optimal solution despite the fact that it was not being measured.

\label{sec:numerical lllustrations}

\section{Conclusions}\label{sec:conclusions}

This paper introduces a framework of nonlinear controllers whose purpose is to drive a given LTI system to the optimal solution of some predefined optimization problem. The controllers are composed of an \textit{estimator}, \textit{optimizer}, and \textit{driver}. We give specific design requirements and possible design choices for each of these modules. Our analysis shows that under mild assumptions, optimal steady-state performance can be guaranteed. Moreover, we give a sufficient condition, in terms of a LMI, such that global exponential asymptotic stability of the optimal steady-state can be guaranteed. Lastly, we present numerical illustrations that demonstrate how the design choices relate to the rate of exponential convergence. The main focus of future work includes further generalizing the proposed framework. In particular, we will consider driving a LTI system to the optimal solution of a \textit{constrained} optimization problem and when there are multiple LTI systems that occur in a distributed setting.

\bibliographystyle{IEEEtran}  
\bibliography{IEEEabrv,references} 

\end{document}